# POO-LPSP: Parallel Osprey Optimized Least Penalty-Squared Prioritization Methods for Priority Derivation in the Analytic Hierarchy Process

**Kevin Kam Fung YUEN**[1*]

[1] School of Science, Monash University Malaysia, Sunway, Malaysia

[*]E-mail: kevin.yuen@monash.edu; kevinkf.yuen@gmail.com;

**Abstract.** Pairwise comparison (PC) via pairwise reciprocal matrices (PRMs) is central to the Analytic Hierarchy Process (AHP). Although the traditional eigenvector method is widely applied to derive priorities, its theoretical robustness in reflecting true priority vectors remains debated. Building upon a previous iteration of this study, this research develops the revised Least Penalty-Squared Prioritization (LPSP) optimization models, including the revised Least Product of Penalty and Direct Squares (LPPDS) and revised Weighted Squares (LPPWS), to minimize the revised Root Mean Penalty-Squared Variance (RMPSV) and the revised Root Mean Penalty-Weighted Square Variance (RMPSWV). However, solving these non-linear formulations is computationally complex for decision-makers. To overcome these limitations, this study proposes the Parallel Osprey Optimized Least Penalty-Squared Prioritization (POO-LPSP) method. By integrating an improved bio-inspired metaheuristic Parallel Osprey Optimization Algorithm (POOA), this framework efficiently solves complex LPSP models to minimize RMPSV and RMPSWV, thereby enhancing prioritization reliability. The practical utility and computational efficiency of the POO-LPSP method are validated through a numerical application focusing on a Generative AI (GAI) vendor selection problem. To extend, POO-LPSP can serve as a robust alternative to Saaty's Eigen system method for AHP applications.

## 1. Introduction

The origins of the Pairwise Comparison (PC) method may be traced back to the pioneering 13th-century social choice and voting theories of the philosopher R. Llull [7]. In 1927, L. L. Thurstone introduced "The Law of Comparative Judgment" [23], establishing a formal psychometric continuum that mathematically converted subjective binary choices into scalable interval measurements. Thomas L. Saaty introduced a structured ratio scale framework in pairwise reciprocal matrices [21] leading to the development of the Analytic Hierarchy Process (AHP) [22], providing a mathematically rigorous foundation for multi-criteria decision-making that uses subjective human perception for priority vector derivation. However, given the existence of diverse PC methodologies, the AHP should not be conflated with the broader concept of PCs [13].

The AHP systematically decomposes a complex decision problem through four core stages: structural definition, assessment, local prioritization, and global synthesis [29]. During the definition stage, decision-makers define the objective, a set of alternatives $T = \{t_1, \dots, t_j, \dots, t_m\}$, and a set of evaluation criteria $C = \{c_1, \dots, c_i, \dots, c_n\}$. In the assessment stage, pairwise comparisons are executed to construct Pairwise Reciprocal Matrices (PRMs). Each PRM is denoted by $A = [a_{ij}]$ such that

$$0 < a_{ij} = a_{ji}^{-1}, \forall i, j \in 1,2, \dots, n \tag{1}$$

The entries of $A$ are assigned values via a 1-to-9 assessment system, measuring the dominance of item $i$ against item $j$. To evaluate PRM consistency, the Consistency Ratio ($CR$) is determined by quotient of the Consistency Index (CI) and Random Index (RI).

$$CR = \frac{CI}{RI} \tag{2}$$

The $CI$ is a function of the matrix dimension $n$ and its principal eigenvalue $\lambda_{max}$, defined as:

$$CI = \frac{\lambda_{max} - n}{n - 1} \tag{3}$$

Furthermore, $\lambda_{max}$ can be calculated by:

$$\lambda_{max} = \frac{1}{n} \sum_{1 \le i < j \le n} \frac{\delta_{ij}^2}{1 + \delta_{ij}}, \delta_{ij} = \left(\frac{a_{ij}}{w_i / w_j} - 1\right) \tag{4}$$

A perfectly consistent PRM always has $\lambda_{max} = n$. The condition $\lambda_{max} < n$ is mathematically impossible, meaning that any deviation from perfect consistency will strictly result in $\lambda_{max} > n$.

In the prioritization stage, extracting a normalized priority vector $W = [w_1, \dots, w_n]$ such that $\sum_{i=1}^{n} w_i = 1$ from a PRM requires applying a Prioritization Operator (PO). The conventional Analytic Hierarchy Process (AHP) utilizes the Eigen System to find $W$ by solving the principal eigenvalue limit problem:

$$w' = \lim_{k \to \infty} \left(\frac{A^k e^T}{e A^k e^T}\right) \quad \Longrightarrow \quad w_i = \frac{w'_i}{\sum_{i=1}^{n} w'_i}, i = 1,2, \dots, n. \tag{5}$$

During the synthesis phase, a weighted average operator aggregates these weighted local priories into a comprehensive global priority vector, $T = [t_1, t_2, \dots, t_j, \dots, t_m]$, which ultimately determines the final ranks for the alternatives.

Despite the wide adoption of the AHP, these four stages frequently introduce operational anomalies, the phenomenon of rank reversal, where the relative preference ordering of alternatives shifts unexpectedly. [29] indicated that these vulnerabilities fundamentally stem from four problem domains: i) the structural selection of criteria in the definition stage; ii) the selection of numerical scales during assessment; iii) the selection of prioritization operators (POs) during prioritization; iv) selection of aggregation operators during final synthesis.

This research specifically investigates Problem iii). The choice of PO is a critical pivot point in multi-criteria decision-making, as different mathematical formulations can yield completely contradictory priority rankings from the exact same underlying PRMs. Although there are various POs ranging from the algebraic forms to optimization forms such as [21; 22],[10], [25], [14], [6], [8], [16], [18], [26], [30] [32], [24], extensive comparative analyses [11] , [19], [5], [29], [12], [31] have established that no single prioritization technique achieves universal superiority when handling inconsistent judgments.

Ultimately, if a unique or novel operator is proposed, evidence must rigorously demonstrate its superior performance against existing methods based on convincing and reasonable measurement criteria. To achieve this, the present research addresses problem 3 by improving a least penalty square prioritization method originally proposed by [28; 30]. As a closed-form solution cannot be analytically derived for this method, numerical optimization is required. The Osprey Optimization Algorithm (OOA) is revised and tailor-made to solve this optimization problem.

OOA represents a metaheuristic approach inspired by osprey predatory behaviors [9]. The algorithm is claimed to mathematically simulate the osprey's dual-phase foraging strategy, which consists of detecting and capturing prey, followed by relocating to consume it. By effectively simulating this behavior to balance global exploration and local exploitation, OOA demonstrates superior optimization performance across various benchmark functions compared to classical well-known algorithms [9]. Recent literature highlights the expanding utility and algorithmic evolution of OOA across diverse engineering and optimization domains. The baseline algorithm has been successfully deployed for lithium-ion battery parameter estimation [1], complex distributed generation placement [20], random forest regression tuning [15] , and joint antenna array vector configurations [27]. To further enhance its mathematical framework, researchers have introduced advanced structural variants, including a Pareto-based OOA (POOA) for multi-output support vector regression[2] and an exploitation-enhanced levy osprey optimization algorithm (LOOA) [17] .

The contribution of the article is presented in several key areas. Section 2 introduces a new PO measurement metric, the Root Mean Penalty-Squared Variance (RMPSV), which integrates violation with squared variance. Section 3 formulates the Least Penalty-Squared Prioritization (LPSP) optimization model to derive an optimized priority vector by minimizing the RMPSV. Section 4 designs a customized Parallel Osprey Optimization Algorithm (POOA) to solve the LPSO models and classical optimization models such as Direct Least Squares (DLS) [6] and Weighted Least Squares (WLS) [6] to calculate priorities (or weights) from pairwise reciprocal matrices. Section 5 validates the practical applicability of the proposed framework via a Generative AI (GAI) service vendor selection case study, benchmarking its performance against established approaches including the Eigenvector, Geometric Mean (GM) [8], DLS, and the Pseudo Inverse Gram Matrix (PIGM) method [31] . Section 6 offers a comprehensive summary of methodological contributions alongside prospective future work.

## 2. Root Mean Penalty-Squared Variance

### *2.1 Foundational Deviation and Variance Metrics*

Selecting the most appropriate Prioritization Operator (PO) requires a robust measurement model supported by empirical validation. PO measurement models evaluate operator fitness, allowing decision-makers to systematically benchmark and select the optimal PO. This section reviews several foundational measurement models and, by synthesizing their advantages while mitigating their limitations, introduces a novel variance framework.

[11] introduced the Total Deviation (TD) to quantify the sum of squared deviations between the derived weight ratios and their corresponding matrix entries. Building upon this, [19] defined Euclidean Distance (ED) as the square root of TD:

$$ED(A, W) = \sqrt{\sum_{i=1}^{n} \sum_{j=1}^{n} \left(a_{ij} - \frac{w_i}{w_j}\right)^2} \tag{6}$$

To enhance computational efficiency via matrix operations, this study reformulates ED using the squared Frobenius norm:

$$ED(A,W) = \| A - \mathrm{w}(\mathrm{w}^{-1})^T \|_F = \sqrt{\sum (A - W(W^{-1})^T)^{\circ 2}} \tag{7}$$

As ED scales with the $n \times n$ dimension of the reciprocal matrix $A$, averaging the values yields a more interpretable metric. [29] proposed the Root Mean Square Variance which calculates the root of the averaged sum of squared deviations:

$$RMSV(A,W) = \sqrt{\frac{1}{n \times n} \sum_{i=1}^{n} \sum_{j=1}^{n} \left( a_{ij} - \frac{w_i}{w_j} \right)^2} \tag{8}$$

To facilitate efficient matrix computation, this paper expresses RMSV as:

$$RMSV(A,W) = \frac{1}{n} \sqrt{\sum (A - W(W^{-1})^T)^{\circ 2}} \tag{9}$$

A variation of this metric is the Root Mean Square-Weighted Variance (RMSWV), defined as:

$$RMSWV(A,W) = \sqrt{\frac{1}{n \times n} \sum_{i=1}^{n} \sum_{j=1}^{n} \left( w_i - a_{ij} w_j \right)^2} \tag{10}$$

In matrix form, RMSWV is calculated as:

$$\begin{aligned} RMSWV(A,W) &= \frac{1}{n} \| W e_n^T - A \times \mathrm{diag}(W) \|_F \\ &= \frac{1}{n} \sqrt{\sum (W e_n^T - A \times \mathrm{diag}(W))^{\circ 2}} \end{aligned} \tag{11}$$

*2.2 Penalty Weights*

A critical limitation of TD, ED, RMSV, and RMSWV is their failure to account for consistency violation penalties. For example, the penalty condition$\left( w_i > w_j \wedge a_{ij} < 1 \wedge a_{ij} \neq \frac{w_i}{w_j} \right)$ is fundamentally different from the condition $\left( w_i > w_j \wedge a_{ij} > 1 \wedge a_{ij} \neq \frac{w_i}{w_j} \right)$. These distinct violations should not be weighed equally. To quantify the variance linked to specific penalty weights, [11] introduced Minimum Violation (MV) index.

$$MV(A,W) = \sum_i \sum_j I_{ij}\,,$$

$$where\ I_{ij} = \begin{cases} 1 & , w_i > w_j \ \&\ a_{ij} < 1 \\ 0.5 & , w_i = w_j \ \&\ a_{ij} \neq 1 \\ 0.5 & , w_i \neq w_j \ \&\ a_{ij} = 1 \\ 0 & , otherwise \end{cases} \tag{12}$$

[29; 30] identified a logical flaw in the original MV piecewise function above, correcting it to ensure $I_{ij}$ equals 1 if $(w_i < w_j \wedge a_{ij} > 1)$. Furthermore, because the raw MV value is dependent on the matrix size ($n^2$), adopting a mean value provides a more stable metric for comparing POs. Based on these improvements, Mean MV (MMV) was formulated. Leveraging reciprocal properties, this study refines the flipped inequality signs of the MMV to properly scale severity:

$$MMV(A,W) = \frac{1}{n^2} \left( \sum_i \sum_j I_{ij} \right) \tag{13}$$

$$\text{Where } I_{ij} = \begin{cases} 1, & (w_i > w_j \wedge a_{ij} < 1) \\ 1, & (w_i < w_j \wedge a_{ij} > 1) \\ 0.5, & (w_i = w_j \wedge a_{ij} \neq 1) \\ 0.5, & (w_i \neq w_j \wedge a_{ij} = 1) \\ 0, & \text{otherwise} \end{cases}$$

*2.3 The revised Framework: Root Mean Penalty-Squared Variance*

A structural pitfall of MMV is that it relies strictly on discrete penalty scores, completely ignoring the magnitude of the variance values. To integrate RMSV/RMSWV variance tracking with MMV structural penalty logic, this study improves upon the hybrid forms from [28; 30] by proposing improved Root Mean Penalty-Squared Variance (RMPSV). This metric is formalized in two variations.

The standard RMPSV *is* proposed as follows.

$$RMPSV(A, W) = \frac{1}{n}\sqrt{\sum P(A - W(W^{-1})^T)^{\circ 2}} = \frac{1}{n}\sqrt{\sum_i \sum_j P_{ij}\left(a_{ij} - \frac{w_i}{w_j}\right)^2} \quad (14)$$

where the penalty matrix, $P = [P_{ij}]$, is determined by:

$$P_{ij} = \begin{cases} \beta_1, & (w_i > w_j \wedge a_{ij} > 1) \\ \beta_1, & (w_i < w_j \wedge a_{ij} < 1) \\ \beta_1, & (w_i = w_j \wedge a_{ij} = 1) \\ \beta_2, & (w_i = w_j \wedge a_{ij} \neq 1) \\ \beta_2, & (w_i \neq w_j \wedge a_{ij} = 1) \\ \beta_3, & \text{otherwise} \end{cases} \quad (15)$$

The weighted variant, Root Mean Penalty-Weighted Square Variance (RMPSWV), is formulated as:

$$\begin{aligned} RMPSWV(A, W) &= \frac{1}{n}\sqrt{\sum P \times (We_n^T - A \times \mathrm{diag}(W))^{\circ 2}} \\ &= \frac{1}{n}\sqrt{\sum_i \sum_j P_{ij}\left(w_i - a_{ij}w_j\right)^2} \end{aligned} \quad (16)$$

Crucially, the baseline consistency condition $(w_i = w_j \wedge a_{ij} = 1)$ for $P_{ij}$, which was omitted in my previous versions [28; 30], is formally integrated here. The penalty weight vector is defined as $\beta = \{\beta_1, \beta_2, \beta_3\}$, subject to the constraint $1 \leq \beta_1 \leq \beta_2 \leq \beta_3$. Under this generalized framework, RMSV and RMSWV are simply special cases of RMPSV and RMPSWV where $\beta_1 = \beta_2 = \beta_3 = 1$. In the historical MMV framework, the implied weights were $\beta_1 = 0, \beta_2 = 0.5, \beta_3 = 1$. However, assigning a zero value to $\beta_1$ inadvertently cancels out valid variance data. To resolve this and ensure mathematical stability, the default penalty settings proposed for this updated model are defined as $\beta_1 = 1, \beta_2 = 3$, and $\beta_3 = 10$.

## 3. Least Penalty-Squared Optimization Prioritization Operators

[6] introduced the Direct Least Squares (DLS) method, which derives priority weights by minimizing the total variance between the reciprocal matrix entry $a_{ij}$ and the weight ratio $\frac{w_i}{w_j}$:

$$\begin{aligned} &Min \quad DS=\sum_{i=1}^{n}\sum_{j=1}^{n}\left(a_{ij}-\frac{w_i}{w_j}\right)^2 \\ &S.T \quad \sum_{i=1}^{n} w_i = 1,\ w_i > 0, i = 1,2,\dots,n . \end{aligned} \tag{17}$$

However, DLS method lacks a closed-form analytical solution, produces alternative optima, which means that different $W$ values leading to the exact same minimized DS value, and is very challenging to resolve via numerical optimization(some discussion can be found in [6][3][4][31]). Thus, [6] transformed the DLS objective function to the Weighted Least Squares (WLS):

$$\begin{aligned} &Min \quad WS=\sum_{i=1}^{n}\sum_{j=1}^{n}\left(w_i - a_{ij}w_j\right)^2 \\ &S.T. \quad \sum_{i=1}^{n} w_i = 1,\ w_i > 0, i = 1,2,\dots,n . \end{aligned} \tag{18}$$

A closed-form solution for WLS remained absent in the literature until [31] introduced four variants of the Inverse Gram Matrix methods. Among these, the Pseudo Inverse Gram Matrix (PIGM) method [31] is formulated as follows:

$$w = \frac{G^{-1}e}{e^T G^{-1} e}, \qquad g_{ij} = \begin{cases} (n-1) + \sum_k a_{kj}^2 & i = j \\ 1 - a_{ij} - a_{ji} & i \neq j \end{cases}, \forall g_{ij} \in G. \tag{19}$$

where $e = (1,\cdots,1)^T$ is a column vector of ones with a size equal to the row dimension of $G^{-1}$.

[8] proposed Logarithm Least Square (LLS) PO has the following form:

$$\begin{aligned} &Min \quad \sum_{i=1}^{n}\sum_{j>i}^{n}\left(ln\, a_{ij} - (ln\, w'_i - ln\, w'_j)\right)^2 \\ &S.T. \quad \prod_{i=1}^{n} w'_i = 1,\ w'_i > 0, i = 1,2,\dots,n \end{aligned} \tag{20}$$

Geometric Mean (GM) is the LLS closed form solution, which produces unnormalized solution $\{w'_i\}$, and the normalization [29; 30] is needed.

$$w'_i = \prod_{j=1}^{n} a_{ij}^{1/n} \implies w_i = \frac{w'_i}{\sum_{i=1}^{n} w'_i}, \quad i = 1,2,\dots,n \tag{21}$$

Least Penalty-Squared Optimization (LPSO) POs mathematically integrate penalty parameters directly into the optimization model to enforce ordinal consistency. Incorporating the updated penalty structures $P_{ij}$, this study refines two primary LPSO operators established in prior iteration of the work [28; 30], namely the Least Product of Penalty and Direct Squares (LPPDS) and the Least Product of Penalty and Weighted Squares (LPPWS). The LPPDS maps the dynamic penalty weights ($P_{ij}$) onto the DLS objective function, yielding the following formulation:

$$\begin{aligned} &Min \quad PPDS=\sum_{i=1}^{n}\sum_{j=1}^{n} P_{ij}\cdot\left(a_{ij}-\frac{w_i}{w_j}\right)^2 \\ &S.T. \quad \sum_{i=1}^{n} w_i = 1,\ w_i > 0, i = 1,2,\dots,n \end{aligned} \tag{22}$$

$P_{ij}$ is shown in Eq. (15). The penalty weights $\beta_1, \beta_2$, and $\beta_3$, by default, are defined as 1,3, and 10 respectively to penalize consistency violations. Similarly, the LPPWS applies the penalty matrix $P_{ij}$ to the WLS formulation, defined as follows:

$$\begin{aligned} &Min \quad \text{PPWS}=\sum_{i=1}^{n}\sum_{j=1}^{n} P_{ij}\cdot\left(w_i - a_{ij}w_j\right)^2 \\ &\text{S.T} \quad \sum_{i=1}^{n} w_i = 1,\ w_i > 0, i = 1,2,\dots,n \end{aligned} \tag{23}$$

Unlike WLS, neither LPPDS nor LPPWS possesses a tractable closed-form solution. Because both models are non-convex, non-differentiable, and governed by piecewise step-functions based on the decision variables, traditional gradient-based solvers, such as Gradient descent, Newton-Raphson or interior-point methods, are fundamentally ineffective. These classical approaches depend on smooth derivatives and will inevitably become trapped in local minima. Consequently, employing a metaheuristic approach is an effective way to reliably explore the global search space. However, calculating the dynamic penalty conditions of $P_{ij}$ for every candidate solution within a metaheuristic population across continuous iterations imposes a severe computational bottleneck.

To overcome these challenges, this study proposes a Parallel Osprey Optimization Algorithm (POOA), inheriting the OOA metaheuristic search mechanism, enable it to successfully escape local minima and approximate the global optimum within a non-convex space. Furthermore, integrating parallel computing architecture allows the algorithm to simultaneously evaluate the complex penalty matrices across the entire candidate population. This parallelization drastically reduces execution time, rendering the LPPDS and LPPWS models highly viable for integration into real-time decision-making frameworks.

## 4. Parallel Osprey Optimization Algorithm

To address the computational limitations and strict structural constraints inherent in decision matrix prioritization, this study proposes Parallel OOA (POOA), an enhanced variant of the Osprey Optimization Algorithm (OOA) [9]. Although the baseline OOA features strong global exploration capabilities, its original design is ill-equipped for this domain due to two primary limitations:

- **Unconstrained Mechanics:** The original algorithm lacks explicit constraint-handling, making it unsuited for the strict mathematical boundaries $(\sum \mathrm{w_i} = 1, \mathrm{w_i} > 0)$ governing priority weight vectors.
- **Sequential Framework:** Its single-threaded, linear execution becomes computationally inefficient when navigating highly non-convex and non-smooth optimization problem.

To achieve real-time operational capability, the proposed POOA significantly advances the standard metaheuristic foundation through four core architectural modifications. The complete POOA execution flow, implemented in R, is detailed in Algorithm 1.

The first major enhancement introduces a dynamic parallel evaluation mechanism to mitigate the severe computational bottlenecks caused by search space to solve the complex objective functions. This multi-core processing layer evaluates population fitness concurrently, governed by an intelligent auto-gatekeeper (Step 2) that dynamically reverts to sequential processing when search dimensions are small, specifically $(N \times m \times S) < 1000$, to prevent operating system thread-spawning overhead. Second, to address the tendency of traditional R-based parallel clusters to inadvertently clone heavy data structures and risk memory leakage, the POOA incorporates a memory-safe environment export filter. Operating within Step 2, this parsing routine utilizes custom scoping functions to strip raw datasets from the active workspace, successfully exporting only core functional components prior to mapping.

Third, the POOA abandons the static maximum iteration limits of classical OOA, which often result in thousands of wasted computational cycles. Instead, it integrates an adaptive early-stopping criterion (Step 12) featuring a stagnation monitor. An active loop-break triggers immediately if a non-improving metric, $no_improvement_counter$, breaches a user-defined threshold $S$. Finally, because priorities satisfy normalization axiom (e.g., summing to 1), the POOA uses an explicit constraint repair mechanism ($f_{repair}$). Implemented during vector initialization (Step 3) and within the core exploration and exploitation loops (Steps 6 and 9), this operator guarantees all candidate positions remain within a mathematically valid.

Algorithm 1: Parallel Osprey Optimization Algorithm (POOA)

| Step | Algorithmic Logic |
|---|---|
| Input: | {Objective function $obj_fn$; search space boundary vectors $lb, ub \in \mathbb{R}^m$; maximum iterations $max_iter$; population size $N$; early-stopping stagnation threshold $S$; hyperplane projection rule $repair_fn$; requested computing cores $n_cores$} |
| Output: | {Normalized optimal priority vector $w^*$, global minimum fitness score $F_{best}$, convergence track.}. |
| Initialization | |
| Step 1: | Define problem dimension: $m \leftarrow length(lb)$. |
| Step 2: | Parallel Infrastructure Configuration Gatekeeper:<br>if $n_cores > 1$, then<br>if $(N \times m \times S) < 1000$, then $n_cores \leftarrow 1$<br>else<br>• Initialize parallel compute cluster: $cl \leftarrow makeCluster(n_cores)$.<br>• Register an execution panic-hook to prevent cluster hangs: $on.exit(stopCluster(cl))$.<br>• Memory-Safe Environment Export Filter: Parse environments (*.GlobalEnv* and *parent.frame()* in R); Identify and strip away heavy data structures or lists; Export only isolated function definitions and scalar rules to worker nodes via *clusterExport* to mitigate system memory leakage.<br>end if<br>end if |
| Step 3: | Population Initialization and Normalization:<br>for each osprey $i = 1$ to $N$ do<br>$X[i,\cdot] \leftarrow lb + rand(1, m) \odot (ub - lb)$.<br>if $repair_fn$ is defined, then<br>Normalize the values in repair function: $X[i,\cdot] \leftarrow repair_fn(X[i,\cdot])$.<br>end for |
| Step 4: | Initial Performance Evaluation:<br>if $n_cores > 1$, then<br>evaluate initial population fitness vector concurrently: $F_vec \leftarrow$ parApply$(cl, X, obj_fn)$. |
| Step 5: | State Tracking Initialization:<br>• Identify absolute best position index: $best_idx \leftarrow \arg\min(F_vec)$.<br>• Extract initial global target: $X_{best} \leftarrow X[best_idx,\cdot]$,<br>• Extract baseline score: $F_{best} \leftarrow F_vec[best_idx]$.<br>• Initialize tracking registers: $no_improvement_counter \leftarrow 0$. |
| Main Loop | Metaheuristic Two-Phase Optimization Loop:<br>for $t = 1$ to $max_iter$ do |
| Step 6: | *PHASE 1*: Position Identification and Hunting (Exploration)<br>for each osprey $i = 1$ to $N$ do<br>Isolate the set of strictly superior candidate solutions:<br>$FP_i \leftarrow \{X[k,\cdot] \mid F_vec[k] < F_vec[i]\}$.<br>if $FP_i$ is empty then<br>Default target set to global leader: $FP_i \leftarrow \{X_{best}\}$.<br>else<br>Append global leader to the candidate set: $FP_i \leftarrow FP_i \cup \{X_{best}\}$.<br>end if<br>Randomly select a target fish position from the set: $SF_i \leftarrow$ sample$(FP_i, 1)$.<br>Generate determination steps $I_{ij} \in \{1,2\}^m$ and stochastic parameters $r_{ij} \in [0,1]^m$.<br>Compute exploration candidate position:<br>$X_{P1}[i,\cdot] \leftarrow X[i,\cdot] + r_{ij} \odot (SF_i - I_{ij} \odot X[i,\cdot])$.<br>Enforce boundary constraints: $X_{P1}[i,\cdot] \leftarrow \max(\min(X_{P1}[i,\cdot], ub), lb)$.<br>if $repair_fn$ is defined, then<br>Map coordinates to priority domain: $X_{P1}[i,\cdot] \leftarrow repair_fn(X_{P1}[i,\cdot])$.<br>end for |
| Step 7: | Phase 1 Evaluation:<br>Concurrently evaluate fitness vector: $F_P1 \leftarrow$ evaluate_population$(X_{P1})$. |
| Step 8: | Phase 1 Greedy Selection Strategy:<br>for each individual osprey $i = 1$ to $N$ do<br>if $F_P1[i] < F_vec[i]$ then update position and fitness:<br>$X[i,\cdot] \leftarrow X_{P1}[i,\cdot]$;<br>$F_vec[i] \leftarrow F_P1[i]$. |
| Step 9: | PHASE 2 (Exploitation): Carrying Fish to Suitable Positions<br>for each individual osprey $i = 1$ to $N$ do |

- Generate localized stochastic sampling array: $r \in [0,1]^m$.
- Compute exploitation position using a time-decaying contraction factor:

  $X_{P2}[i,\cdot] \leftarrow X[i,\cdot] + \frac{lb + r \odot (ub - lb)}{t}$.
- Enforce boundary constraints:

  $X_{P2}[i,\cdot] \leftarrow \max(\min(X_{P2}[i,\cdot], ub), lb)$.
- if $repair_fn$ is defined, then $X_{P2}[i,\cdot] \leftarrow repair_fn(X_{P2}[i,\cdot])$.

end for

Step 10: Phase 2 Greedy Evaluation

Concurrently Phase 2 candidate fitness array:

$F_P2 \leftarrow$ evaluate_population($X_{P2}$).

Step 11: Phase 2 Greedy Selection Strategy:

for each individual osprey $i = 1$ to $N$ do

if $F_P2[i] < F_vec[i]$, then update position and fitness:

$X[i,\cdot] \leftarrow X_{P2}[i,\cdot]$;

$F_vec[i] \leftarrow F_P2[i]$.

end for

Step 12: Global Memory Management & Adaptive Stopping Check:

*current_best_idx←argmin(F_vec);*

$current_best_F \leftarrow F_vec[current_best_idx]$;

if $current_best_F < F_{best}$ then

Record new global optimum:

- $F_{best} \leftarrow current_best_F$;
- $X_{best} \leftarrow X[current_best_idx,\cdot]$;

Reset stagnation tracking register: $no_improvement_counter \leftarrow 0$.

else

Increment stagnation tracking register:

$no_improvement_counter \leftarrow no_improvement_counter + 1$.

end if

if $no_improvement_counter \geq S$, then

Terminate optimization loop prematurely due to Stagnation Limit Reached.

end if

end for

Step 13: Normalization & Cleanup:

If Cluster instance $cl$ is active, then execute $stopCluster(cl)$.

Map the unbounded optimal array into an exactly normalized priority weight vector:

$$w^* \leftarrow \frac{X_{best}}{\sum_{k=1}^{m} X_{best,k}}$$

return $\{w^*, F_{best}$, convergence track$\}$

## 5. Application of Generative AI Service Provider Selection

A prominent healthcare enterprise leverages Generative AI (GAI) to automate administrative workflows, clinical documentation, and unstructured literature synthesis. This minimizes clinicians' cognitive burdens, accelerating research and diagnostics while maintaining data security. Choosing an optimal vendor is vital to balance clinical efficacy, technological agility, and operational efficiency. To address this vendor selection challenge, this study utilizes the enhanced AHP with seven POs, Eigenvector, Geometric Mean (a closed form of LLM), DLS, PIGM (a closed form of WLS), LPPWS and LPPDS, proposed in Section 3, with respect to the five PO metrics, MMV, RMSV, RMSWV, RMPSV and RMPSWV, proposed in Section 4.

*5.1 Decision Problem Structure*

The hierarchical decision model, presented in Figure 1, is designed to identify an optimal GAI service provider for sustainable business operations. The evaluation framework comprises six core criteria ($C_1$–$C_6$) applied to a shortlist of five anonymized alternatives ($t_1$–$t_5$). For the purposes of vendor neutrality and data masking, these alternatives are representative of major platforms (e.g., Google's Gemini 2.0, OpenAI's GPT-5.5, Claude 5, Meta's Llama 4, Cohere's Command A+, DeepSeek's DeepSeek-V4, xAI's Grok 4.5, Mistral Large 3, etc.). The operational definitions of these criteria are provided below.

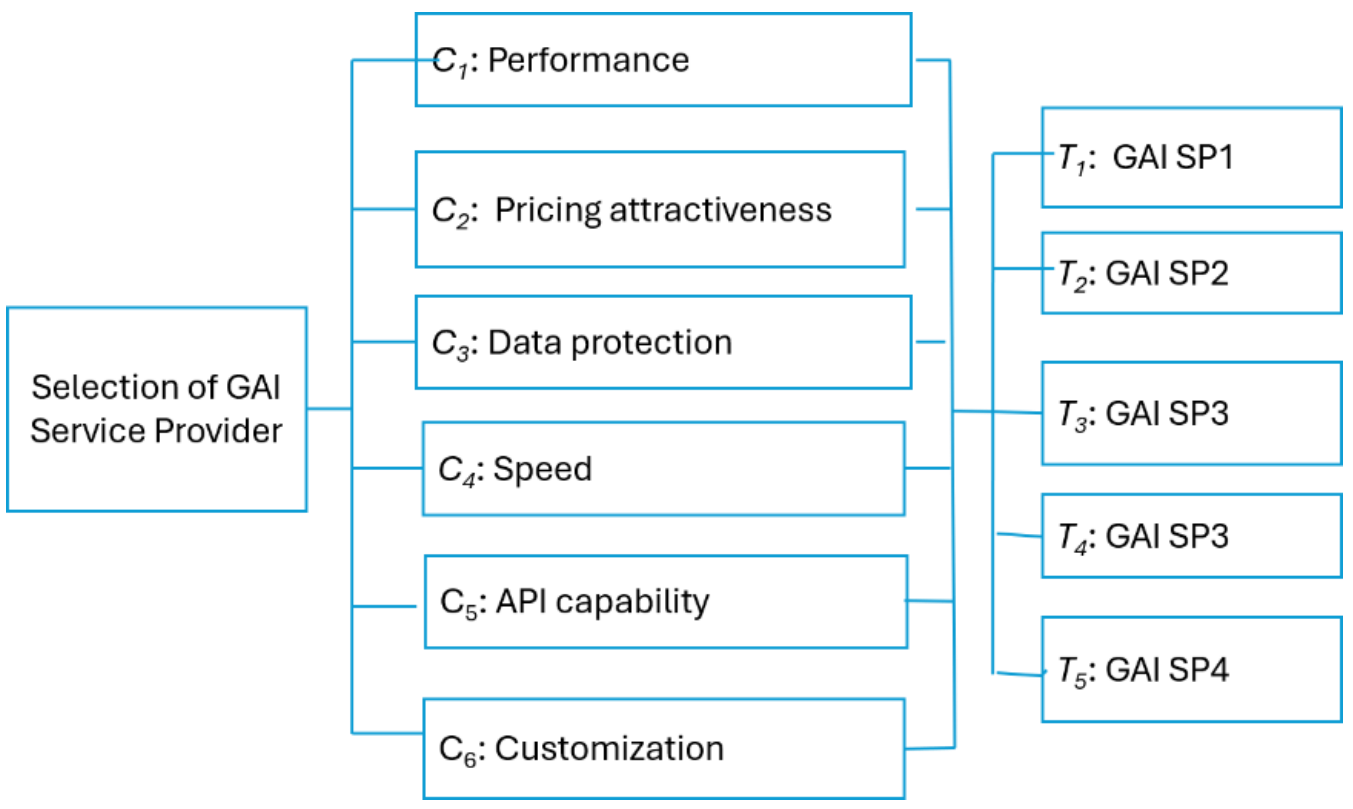

**Figure 1.** Structure for selecting the best GAI service provider

- Performance ($C_1$): evaluates cognitive reasoning in complex biomedical tasks, measuring hallucination rates, semantic coherence, and diagnostic accuracy to ensure clinical fidelity.
- Pricing attractiveness ($C_2$): assesses total cost of ownership across token pricing schedules, subscription models, enterprise volume discounts, and batch-processing cost-mitigation structures.
- Data protection ($C_3$): verifies privacy and security compliance for protected health information (PHI), encompassing enterprise encryption, data residency, zero-retention training policies, and cybersecurity resilience.

- Speed ($C_4$): quantifies real-time responsiveness via time-to-first-token (TTFT), throughput velocity, generation speed, and SLA-backed API uptime guarantees.
- API capability ($C_5$): gauges interoperability via developer documentation, SDK robustness, and integration compatibility with legacy enterprise infrastructure and clinical software tools.
- Customization ( $C_6$ ): evaluates a system's adaptability to clinical terminology through retrieval-augmented generation (RAG), hyperparameter tuning, and parameter-efficient fine-tuning (PEFT).

**Table 1.** PRM for Six Criteria ($A_0$ , CR= 0.094)

| $A_0$ | $C_1$ | $C_2$ | $C_3$ | $C_4$ | $C_5$ | $C_6$ |
|---|---|---|---|---|---|---|
| $C_1$ | 1 | 2 | 1/4 | 6 | 5 | 2 |
| $C_2$ | 1/2 | 1 | 1/6 | 3 | 1/3 | 1/4 |
| $C_3$ | 4 | 6 | 1 | 8 | 7 | 5 |
| $C_4$ | 1/6 | 1/3 | 1/8 | 1 | 1/4 | 1/5 |
| $C_5$ | 1/5 | 3 | 1/7 | 4 | 1 | 1/3 |
| $C_6$ | 1/2 | 4 | 1/5 | 5 | 3 | 1 |

**Table 2.** PRMs for priorities among five alternatives with respect to each criterion.

| | $A_1$ (CR= 0.019) | | | | | $A_2$ (CR = 0.020) | | | | |
|---|---|---|---|---|---|---|---|---|---|---|
| | $T_1$ | $T_2$ | $T_3$ | $T_4$ | $T_5$ | $T_1$ | $T_2$ | $T_3$ | $T_4$ | $T_5$ |
| $T_1$ | 1 | 2 | 1 | 4 | 5 | 1 | 1/3 | 1/3 | 1/5 | 1/6 |
| $T_2$ | 1/2 | 1 | 1/2 | 2 | 3 | 3 | 1 | 2 | 1/2 | 1/3 |
| $T_3$ | 1 | 2 | 1 | 2 | 3 | 3 | 1/2 | 1 | 1/3 | 1/4 |
| $T_4$ | 1/4 | 1/2 | 1/2 | 1 | 2 | 5 | 2 | 3 | 1 | 1/2 |
| $T_5$ | 1/5 | 1/3 | 1/3 | 1/2 | 1 | 6 | 3 | 4 | 2 | 1 |
| | $A_3$ (CR=0.017) | | | | | $A_4$ (CR=0.024) | | | | |
| $T_1$ | 1 | 1 | 1 | 4 | 5 | 1 | 1 | 2 | 1/4 | 1/5 |
| $T_2$ | 1 | 1 | 1 | 2 | 5 | 1 | 1 | 3 | 1/2 | 1/3 |
| $T_3$ | 1 | 1 | 1 | 3 | 5 | 1/2 | 1/3 | 1 | 1/4 | 1/5 |
| $T_4$ | 1/4 | 1/2 | 1/3 | 1 | 3 | 4 | 2 | 4 | 1 | 1 |
| $T_5$ | 1/5 | 1/5 | 1/5 | 1/3 | 1 | 5 | 3 | 5 | 1 | 1 |
| | $A_5$ (CR=0.036) | | | | | $A_6$ (CR = 0.051) | | | | |
| $T_1$ | 1 | 3 | 1/2 | 5 | 5 | 1 | 1/2 | 1/6 | 1/5 | 1/4 |
| $T_2$ | 1/3 | 1 | 1/2 | 3 | 3 | 2 | 1 | 1/2 | 1/4 | 1/3 |
| $T_3$ | 2 | 2 | 1 | 4 | 4 | 6 | 2 | 1 | 1/3 | 1/2 |
| $T_4$ | 1/5 | 1/3 | 1/4 | 1 | 1 | 5 | 4 | 3 | 1 | 3 |
| $T_5$ | 1/5 | 1/3 | 1/4 | 1 | 1 | 4 | 3 | 2 | 1/3 | 1 |

Table 1 delineates the initial PRM formulated to assess the comparative significance of the six criteria. Building upon this, Table 2 provides the specific matrices constructed to rank the five alternatives under each criterion. Given that the Consistency Ratio (CR) across all computed matrices falls safely beneath the guided 0.1 limit, within the acceptance range.

**Table 3.** Weighted Decision matrices and final results with respect to each PO.

| | Eigen | | | | | | | GM/LLS | | | | | | |
|---|---|---|---|---|---|---|---|---|---|---|---|---|---|---|
| | $C_1$ | $C_2$ | $C_3$ | $C_4$ | $C_5$ | $C_6$ | | $C_1$ | $C_2$ | $C_3$ | $C_4$ | $C_5$ | $C_6$ | |
| W | 0.2 | 0.06 | 0.483 | 0.029 | 0.081 | 0.148 | Final | 0.197 | 0.059 | 0.486 | 0.03 | 0.078 | 0.151 | Final |
| $T_1$ | 0.353 | 0.051 | 0.302 | 0.102 | 0.328 | 0.053 | **0.257** | 0.355 | 0.051 | 0.299 | 0.101 | 0.326 | 0.052 | 0.255 |
| $T_2$ | 0.183 | 0.156 | 0.258 | 0.139 | 0.168 | 0.091 | 0.202 | 0.184 | 0.156 | 0.26 | 0.139 | 0.171 | 0.095 | 0.204 |
| $T_3$ | 0.282 | 0.104 | 0.28 | 0.062 | 0.369 | 0.18 | 0.256 | 0.279 | 0.103 | 0.282 | 0.061 | 0.363 | 0.179 | **0.255** |
| $T_4$ | 0.112 | 0.267 | 0.11 | 0.319 | 0.068 | 0.442 | 0.172 | 0.112 | 0.268 | 0.108 | 0.32 | 0.07 | 0.439 | 0.172 |
| $T_5$ | 0.069 | 0.422 | 0.05 | 0.378 | 0.068 | 0.234 | 0.114 | 0.069 | 0.422 | 0.05 | 0.379 | 0.07 | 0.236 | 0.115 |
| | *DLS* | | | | | | | PIGM/WLS | | | | | | |
| W | 0.235 | 0.07 | 0.41 | 0.044 | 0.061 | 0.181 | Final | 0.165 | 0.069 | 0.537 | 0.043 | 0.065 | 0.121 | Final |
| $T_1$ | 0.378 | 0.059 | 0.301 | 0.084 | 0.365 | 0.057 | **0.253** | 0.352 | 0.059 | 0.297 | 0.083 | 0.295 | 0.056 | 0.251 |
| $T_2$ | 0.195 | 0.158 | 0.262 | 0.145 | 0.181 | 0.096 | 0.199 | 0.175 | 0.148 | 0.271 | 0.139 | 0.152 | 0.102 | 0.213 |
| $T_3$ | 0.249 | 0.104 | 0.281 | 0.072 | 0.31 | 0.264 | 0.251 | 0.299 | 0.1 | 0.285 | 0.07 | 0.406 | 0.163 | **0.258** |
| $T_4$ | 0.106 | 0.292 | 0.102 | 0.31 | 0.073 | 0.346 | 0.167 | 0.103 | 0.255 | 0.093 | 0.332 | 0.074 | 0.478 | 0.162 |
| $T_5$ | 0.074 | 0.387 | 0.055 | 0.389 | 0.073 | 0.236 | 0.131 | 0.073 | 0.438 | 0.054 | 0.376 | 0.074 | 0.202 | 0.117 |
| | *LPPDS* | | | | | | | LPPWS | | | | | | |
| W | 0.234 | 0.066 | 0.409 | 0.044 | 0.066 | 0.183 | Final | 0.166 | 0.067 | 0.537 | 0.043 | 0.067 | 0.121 | Final |
| $T_1$ | 0.365 | 0.059 | 0.296 | 0.088 | 0.337 | 0.058 | 0.247 | 0.338 | 0.059 | 0.29 | 0.087 | 0.295 | 0.056 | 0.246 |
| $T_2$ | 0.194 | 0.158 | 0.266 | 0.133 | 0.182 | 0.095 | 0.2 | 0.175 | 0.148 | 0.279 | 0.131 | 0.152 | 0.102 | 0.217 |
| $T_3$ | 0.262 | 0.104 | 0.281 | 0.073 | 0.337 | 0.25 | **0.254** | 0.313 | 0.1 | 0.285 | 0.07 | 0.406 | 0.163 | **0.261** |
| $T_4$ | 0.106 | 0.292 | 0.102 | 0.319 | 0.073 | 0.348 | 0.168 | 0.102 | 0.255 | 0.093 | 0.346 | 0.074 | 0.478 | 0.162 |
| $T_5$ | 0.073 | 0.387 | 0.055 | 0.387 | 0.073 | 0.25 | 0.132 | 0.072 | 0.438 | 0.054 | 0.367 | 0.074 | 0.202 | 0.115 |

### *5.2 Comparisons and discussion*

To resolve the complex, non-linear objective functions underlying DLS, WLS, LPPDS, and LPPWS, the POOA is implemented as the primary metaheuristic solver. POOA efficiently explores the continuous priority weight space to avoid local minima and reliably discover global optimal solutions under varying penalty constraints. Crucially, the POOA framework configured for the WLS objective function yields priority vectors nearly identical to the closed-form PIGM solution, validating the algorithm's high convergence precision and mathematical consistency. Table 3 details the weighted decision matrices and final rankings derived across different POs. For selection outcome of the best provider, GM, PIGM, LPPDS, and LPPWS yield $T_3$ , whereas Eigen and DLS yield $T_1$. That means you may make a wrong best choice if you use classical Eigen method for AHP. Conversely, all evaluated methods consistently rank $T_5$ as the least optimal alternative.

Table 4 reports error measurements, bolding the minimum values for each metric. All PRMs exhibiting a CR less than 0.1 are structurally valid. However, non-zero MMV values indicate that not all POs consistently generate feasible solutions.

**Table 4.** Error measurement metrics.

| $A_0$ (CR= 0.094) | | | | | | $A_1$ (CR= 0.019) | | | | |
|---|---|---|---|---|---|---|---|---|---|---|
| | MMV | RMSV | RMSWV | RMPSV | RMPSWV | MMV | RMSV | RMSWV | RMPSV | RMPSWV |
| Eigen | **0** | 1.685 | 0.100 | 1.685 | 0.100 | 0.04 | 0.345 | 0.042 | 0.357 | 0.050 |
| GM | **0** | 1.665 | 0.098 | 1.665 | 0.098 | 0.04 | 0.338 | 0.042 | 0.352 | 0.052 |
| PIGM | 0.056 | 1.304 | **0.075** | 1.699 | 0.106 | 0.04 | 0.371 | **0.038** | 0.377 | 0.044 |
| DLS | 0.056 | **0.998** | 0.134 | 1.510 | 0.155 | 0.04 | **0.281** | 0.051 | 0.331 | 0.073 |
| LPPWS | **0** | 1.306 | 0.075 | 1.306 | **0.075** | 0.04 | 0.419 | 0.040 | 0.420 | **0.041** |
| LPPDS | **0** | 1.005 | 0.135 | **1.005** | 0.135 | 0.04 | 0.290 | 0.046 | **0.321** | 0.062 |
| $A_2$ (CR = 0.020) | | | | | | $A_3$ (CR=0.017) | | | | |
| Eigen | 0 | 0.519 | 0.041 | 0.519 | 0.041 | 0.12 | 0.402 | 0.038 | 0.409 | 0.044 |
| GM | 0 | 0.531 | 0.042 | 0.531 | 0.042 | 0.12 | 0.400 | 0.037 | 0.406 | 0.041 |
| PIGM | 0 | 0.463 | **0.036** | 0.463 | **0.036** | 0.12 | 0.375 | **0.031** | 0.378 | 0.033 |
| DLS | 0 | **0.359** | 0.053 | **0.359** | 0.053 | 0.12 | **0.358** | 0.033 | 0.365 | 0.038 |
| LPPWS | 0 | 0.463 | **0.036** | 0.463 | **0.036** | 0.12 | 0.386 | **0.031** | 0.387 | **0.032** |
| LPPDS | 0 | **0.359** | 0.053 | **0.359** | 0.053 | 0.12 | 0.359 | 0.032 | **0.363** | 0.036 |
| $A_4$ (CR=0.024) | | | | | | $A_5$ (CR=0.036) | | | | |
| Eigen | 0.080 | 0.500 | 0.046 | 0.521 | 0.054 | **0** | 0.525 | 0.081 | 0.525 | 0.081 |
| GM | 0.080 | 0.499 | 0.045 | 0.521 | 0.054 | **0** | 0.487 | 0.081 | 0.487 | 0.081 |
| PIGM | 0.080 | 0.381 | **0.035** | 0.443 | 0.045 | **0** | 0.644 | **0.073** | 0.644 | **0.073** |
| DLS | 0.080 | **0.360** | 0.038 | 0.442 | 0.055 | 0.120 | **0.376** | 0.102 | 0.885 | 0.300 |
| LPPWS | 0.080 | 0.418 | 0.037 | 0.452 | **0.042** | 0.040 | 0.644 | **0.073** | 0.644 | **0.073** |
| LPPDS | 0.080 | 0.376 | 0.037 | **0.421** | 0.050 | 0.040 | 0.411 | 0.089 | **0.411** | 0.089 |
| $A_6$ (CR = 0.051) | | | | | | Sum of variance | | | | |
| Eigen | 0 | 0.930 | 0.082 | 0.930 | 0.082 | **0.240** | 4.906 | 0.429 | 4.946 | 0.451 |
| GM | 0 | 0.932 | 0.082 | 0.932 | 0.082 | **0.240** | 4.853 | 0.426 | 4.894 | 0.449 |
| PIGM | 0 | 1.005 | **0.072** | 1.005 | **0.072** | 0.296 | 4.545 | **0.359** | 5.010 | 0.408 |
| DLS | 0.080 | **0.672** | 0.141 | 1.016 | 0.242 | 0.496 | **3.406** | 0.552 | 4.908 | 0.915 |
| LPPWS | 0 | 1.005 | **0.072** | 1.005 | **0.072** | 0.280 | 4.644 | 0.364 | 4.679 | **0.370** |
| LPPDS | 0 | 0.678 | 0.135 | **0.678** | 0.135 | 0.280 | 3.479 | 0.529 | **3.558** | 0.560 |

Regarding the specific error metrics shown in Table 4, DLS consistently minimizes RMSV, while PIGM (WLS) minimizes RMSWV. The penalty-based models, LPPDS and LPPWS, consistently produce the lowest RMPSV and RMPSWV, respectively. When MMV equals zero, DLS converges with LPPDS, and PIGM (WLS) converges with LPPWS; they produce identical corresponding error values due to equivalent priority outputs. The overall variance sums confirm this behavior: DLS (3.406) and PIGM (0.359) minimize RMSV and RMSWV respectively, whereas LPPDS (3.558) and LPPWS (0.370) effectively minimize RMPSV and RMPSWV respectively.

In summary, while vendor rankings fluctuate based on the applied PO, $T_3$ emerges as the most robust choice across most mathematical methods. Furthermore, the error analysis validates the efficacy of the proposed LPPDS and LPPWS optimization models. By successfully minimizing penalty-weighted variances while maintaining structural consistency, these models provide a highly reliable prioritization framework for decision-making scenarios where feasibility violations must be strictly penalized.

## 6. Conclusion

While the mathematical robustness of Saaty's eigenvector method for PRMs remains a subject of ongoing academic debate, this study introduces the Parallel Osprey Optimized Least Penalty-Squared Prioritization method to resolve the computational complexity and opacity inherent in optimization-based prioritization within the Analytic Hierarchy Process. By leveraging the bio-inspired POOA, the proposed framework efficiently solves intricate LPSP formulations to derive accurate priority vectors while explicitly minimizing RMPSV and RMPSWV.

A practical enterprise application focusing on Generative AI vendor selection validates the operational feasibility, clarity, and computational efficiency of the POO-LPSP method, establishing it as a highly rigorous and scalable decision-support tool for complex expert evaluations. The experimental results based on application demonstrate that the POOA framework serves as a powerful, versatile metaheuristic solver capable of discovering global optimal solutions across LPSP, WLS, and DLS optimization landscapes. Furthermore, the comparative error analysis indicates that eigenvector method cannot select the best one choice when evaluated against the demonstrated metrics, emphasizing the measurement of penalty and variance.

Future research may focus on deploying the POO-LPSP framework across broader industrial applications, actively replacing classical AHP in high-stakes domains where traditional eigenvector methods yield suboptimal variance resolution. Furthermore, extending the algorithm's highly scalable parallel architecture to support large-scale group decision-making may be explored to efficiently aggregate massive, highly diverse expert evaluation matrices. Finally, integrating POO-LPSP with other MCDM models such as TOPSIS or PROMETHEE may be investigated to further enhance its versatility and cross-disciplinary applicability.

## References

[1] A.H. Alqahtani, H.M. Fahmy, H.M. Hasanien, M. Tostado-Véliz, A. Alkuhayli, and F. Jurado, Parameters estimation and sensitivity analysis of lithium-ion battery model uncertainty based on osprey optimization algorithm, *Energy* **304** (2024), 132204.

[2] M. Alrbai, S. Al-Dahidi, H. Alahmer, L. Al-Ghussain, R. Al-Rbaihat, H. Hayajneh, and A. Alahmer, Integration and Optimization of a Waste Heat Driven Organic Rankine Cycle for Power Generation in Wastewater Treatment Plants, *Energy* **308** (2024), 132829.

[3] S. Bozóki, Solution of the least squares method problem of pairwise comparison matrices, *Central European Journal of Operations Research* **16** (2008), 345-358.

[4] E. Carrizosa and F. Messine, An exact global optimization method for deriving weights from pairwise comparison matrices, *Journal of Global Optimization* **38** (2007), 237-247.

[5] E.U. Choo and W.C. Wedley, A common framework for deriving preference values from pairwise comparison matrices, *Computers & Operations Research* **31** (2004), 893-908.
[6] A.T.W. Chu, R.E. Kalaba, and K. Spingarn, A comparison of two methods for determining the weights of belonging to fuzzy sets, *Journal of Optimization Theory and Applications* **27** (1979), 531-538.
[7] J.M. Colomer, Ramon Llull: from 'Ars electionis' to social choice theory, *Social Choice and Welfare* **40** (2013), 317-328.
[8] G. Crawford and C. Williams, A note on the analysis of subjective judgment matrices, *Journal of Mathematical Psychology* **29** (1985), 387-405.
[9] M. Dehghani and P. Trojovský, Osprey optimization algorithm: A new bio-inspired metaheuristic algorithm for solving engineering optimization problems, *Frontiers in Mechanical Engineering* **Volume 8 - 2022** (2023).
[10] S.I. Gass and T. Rapcsák, Singular value decomposition in AHP, *European Journal of Operational Research* **154** (2004), 573-584.
[11] B. Golany and M. Kress, A multicriteria evaluation of methods for obtaining weights from ratio-scale matrices, *European Journal of Operational Research* **69** (1993), 210-220.
[12] J. Gyani, A. Ahmed, and M.A. Haq, MCDM and Various Prioritization Methods in AHP for CSS: A Comprehensive Review, *IEEE Access* **10** (2022), 33492-33511.
[13] W.W. Koczkodaj, L. Mikhailov, G. Redlarski, M. Soltys, J. Szybowski, G. Tamazian, E. Wajch, and K.K.F. Yuen, Important Facts and Observations about Pairwise Comparisons, *Fundamenta Informaticae* **144** (2016), 1-17.
[14] G. Kou and C. Lin, A cosine maximization method for the priority vector derivation in AHP, *European Journal of Operational Research* **235** (2014), 225-232.
[15] Q. Liang, Y. Xu, B. Xu, and Y. Du, Parameter optimization for in-situ synthesized TiB2/TiC particle composite coatings by laser cladding based on OOA-RFR and U-NSGA-III, *Optics & Laser Technology* **181** (2025), 111755.
[16] C.C. Lin, An enhanced goal programming method for generating priority vectors, *Journal of the Operational Research Society* **57** (2006), 1491-1496.
[17] Manisha and P. Kumar, An improved osprey optimization algorithm to analyse the steady state performance of stainless steel utensil manufacturing unit, *Life Cycle Reliability and Safety Engineering* **15** (2026), 299-316.
[18] L. Mikhailov, A fuzzy programming method for deriving priorities in the analytic hierarchy process, *Journal of the Operational Research Society* **51** (2000), 341-349.
[19] L. Mikhailov and M.G. Singh, Comparison analysis of methods for deriving priorities in the analytic hierarchy process, in: *IEEE SMC'99 Conference Proceedings. 1999 IEEE International Conference on Systems, Man, and Cybernetics (Cat. No.99CH37028)*, 1999, pp. 1037-1042 vol.1031.
[20] P. Rajakumar, P.M. Balasubramaniam, A. Munimathan, M. Alwetaishi, M.I. Gulbarga, and C. Senthilpari, Osprey optimization algorithm for distributed generation integration in a radial distribution system for power loss reduction, *Scientific Reports* **15** (2025), 18775.
[21] T.L. Saaty, A scaling method for priorities in hierarchical structures, *Journal of Mathematical Psychology* **15** (1977), 234-281.
[22] T.L. Saaty, Analytic Hierarchy Process: Planning, Priority, Setting, Resource Allocation, *McGraw-Hill, New York* (1980).
[23] L. Thurstone, A law of comparative judgment, *Psychological Review* **34** (1927), 273-286.
[24] H. Wang, Y. Peng, and G. Kou, A two-stage ranking method to minimize ordinal violation for pairwise comparisons, *Applied Soft Computing* **106** (2021), 107287.
[25] Y.-M. Wang, C. Parkan, and Y. Luo, Priority estimation in the AHP through maximization of correlation coefficient, *Applied Mathematical Modelling* **31** (2007), 2711-2718.
[26] Y.-M. Wang, C. Parkan, and Y. Luo, A linear programming method for generating the most favorable weights from a pairwise comparison matrix, *Computers & Operations Research* **35** (2008), 3918-3930.
[27] H. Xie, C. Wan, D. Quan, G. Gui, and Y. Yao, Movable Antenna-Enhanced Beamforming via an Improved Osprey Optimization Algorithm, *IEEE Wireless Communications Letters* **15** (2026), 3104-3108.
[28] K.K.F. Yuen, Cognitive network process with fuzzy soft computing technique for collective decision aiding, *The Hong Kong Polytechnic University* **Ph.D. thesis** (2009).
[29] K.K.F. Yuen, Analytic hierarchy prioritization process in the AHP application development: A prioritization operator selection approach, *Applied Soft Computing Journal* **10** (2010), 975-989.
[30] K.K.F. Yuen, The Least Penalty Optimization Prioritization Operators for the Analytic Hierarchy Process: A Revised Case of Medical Decision Problem of Organ Transplantation, *Systems Engineering* **17** (2014), 442-461.
[31] K.K.F. Yuen, Inverse gram matrix methods for prioritization in analytic hierarchy process: Explainability of weighted least squares optimization method, *arXiv preprint arXiv:2401.01190* (2024).
[32] J. Zhang, G. Kou, Y. Peng, and Y. Zhang, Estimating priorities from relative deviations in pairwise comparison matrices, *Information Sciences* **552** (2021), 310-327.